\documentclass[twocolumn]{arxiv}

\usepackage{natbib}
\usepackage{amssymb}
\usepackage{amsmath}
\usepackage{graphicx}
\graphicspath{ {fig/} }
\usepackage{epstopdf}

\begin{document}

\begin{frontmatter}

\title{Optimal output consensus for linear systems: A topology free approach\thanksref{footnoteinfo}}

\thanks[footnoteinfo]{The authors gratefully acknowledge the financial support
from the Swedish National Space Technology Research Programme (NRFP)
and the Swedish Research Council (VR).}

\author[Johan]{Johan Thunberg}\ead{johan.thunberg@math.kth.se},
\author[Johan]{Xiaoming Hu}\ead{hu@math.kth.se}

\address[Johan]{KTH Royal Institute of Technology,  S-100 44 Stockholm, Sweden}


\begin{abstract}
In this paper, for any homogeneous system of agents with linear continuous time 
dynamics, we formulate an optimal control problem. In this problem a convex 
cost functional of the control signals of the agents shall be minimized, while
the outputs of the agents shall coincide at some given finite time. This is an 
instance of the rendezvous or finite time consensus problem. We solve this problem 
without any constraints on the communication topology and provide a solution as an explicit feedback control law for the case when the dynamics of the agents is output controllable. It turns out that the communication graph topology induced by the solution is complete. Based on this solution for the finite time consensus problem, we provide a solution to the case of infinite time horizon. Furthermore, we investigate under what circumstances it is possible to express the controller as a feedback control law of the output instead of the states. 
\end{abstract}

\end{frontmatter}

\section{Introduction}
Here we study the output
consensus problem~\cite{kim2011output,xi2012output,xi2012olle} for a homogeneous system of agents with linear dynamics, both in 
finite time and in the asymptotic case (as time tends to infinity).
In the finite time case (rendezvous)~\cite{wang2010finite,sundaram2007finite},
the outputs for all the agents shall be the same 
at some predefined finite time. It is easy to show that
for homogeneous systems of agents with linear dynamics, it is not possible
to construct a linear, time-invariant feedback control law 
based on relative information such that the agents reach consensus in their
states in finite time. With relative information in this context, we 
are referring to pairwise differences between the states of the agents. 

Regarding the output consensus problem, using a decomposition of the state space, we show that
if the dynamics for the agents is output controllable and the nullspace of the 
matrix which maps the state to the output satisfies a certain invariance condition,
there cannot exist a linear continuous in state, time-invariant feedback control law that solves the problem
while using only relative output information in the form of pairwise differences
between the outputs of the agents.
If only relative
information is used, the control laws need to be either
time-varying or non-Lipschitz in order to solve the finite time consensus problem.
Furthermore, the output controllability is a standing assumption in order
to guarantee consensus for arbitrary initial conditions. 

In \cite{cao2009lqr,cao2010optimal}, an optimal linear consensus problem 
is addressed for systems of mobile agents with single-integrator
dynamics. In this setting, the authors constrain the agents to use
only relative information in their controllers, \emph{i.e.},~the
controller of each agent consists of a weighted sum of the differences
between its state and the states of its neighbors. In this setting the
authors show that the graph Laplacian matrix used in the optimal
controller for the system corresponds to a complete directed
graph. There is another line of research on the optimal consensus
problem \cite{semsar2008optimal}, in which the consensus requirement is
reflected in the cost function. However, with such a formulation the
optimal controller in general can not be implemented with relative
state information only.

We formulate the
consensus problem as an optimal control problem, where there are no
restrictions on the controllers besides that the agents shall reach
consensus at some predefined time. 
Note that we do not impose any communication topology on the system
beforehand, instead we are interested in the topology generated by the
solution to the optimal control problem. By solving the problem, we
show that the  
optimal controller uses only relative information. Moreover, the
connectivity graph needs to be completely connected. Thus, for
any other topology between the agents than the complete graph, any
controller constructed will be suboptimal. 
The provided control laws are given in closed form and are bounded and
continuous. The input and output dimensions are arbitrary. 

 Not surprisingly, the
optimal controller requires the measurement of state errors in general. We identify cases where the optimal controller is only
based on the output errors. We also show that in the asymptotic case, there is a corresponding observer based
controller, that is only based on the output errors.

Regarding the theoretical contribution of this work, we use linear vector space optimization 
methods in order to solve the consensus problems. We show that the 
problem can be posed as a certain minimum norm problem in a Hilbert space~\cite{luenberger1997optimization}. 
In this framework the finite time consensus problem is a solution of an optimization problem in this Hilbert space.

\section{Preliminaries}\label{problem}
We consider a system of $N$ agents, where each agent $i$ in the system has the dynamics 
\begin{align*}\label{eq:dyn:problem1}
\dot{x}_i & =  Ax_i(t) + Bu_i(t), \\
y_i & =  Cx_i.
\end{align*}
The variable $x_i(t_0) = x_0$, $x_i(t):\mathbb{R} \rightarrow \mathbb{R}^n$,
$u_i(t):\mathbb{R} \rightarrow \mathbb{R}^m$ and $y_i(t):\mathbb{R}
\rightarrow \mathbb{R}^p$, $A \in \mathbb{R}^{n \times n}$, $B \in
\mathbb{R}^{n \times m}$ and $C \in \mathbb{R}^{p \times n}$. 
It is assumed that $B$
and $C$ are full rank matrices and that the system is output controllable. Let us define 
\begin{align*}
x(t) & = [x_1^T(t), x_2^T(t), \ldots , x_N^T(t)]^T \in \mathbb{R}^{nN}, \\
u(t) & = [u_1^T(t), u_2^T(t), \ldots, u_N^T(t)]^T \in \mathbb{R}^{mN}, \\
y(t) & = [y_1^T(t), y_2^T(t), \ldots, y_N^T(t)]^T \in \mathbb{R}^{pN},
\end{align*}
and the vector
$$a = [a_1, a_2,...,a_N]^T,$$
where $a_i \in \mathbb{R}^{+}$ for all $i$. For each agent $i$, the positive scalar $a_i$ determines how much the control 
signal of agent $i$ shall be penalized in the objective function 
of the optimal control problem; see \eqref{eq:objective} further down.

The matrix 
\begin{equation}
L(a) = -\left(\left(\sum_{i=1}^Na_i\right)^{-1} 1_Na^T - \text{diag}\left([1, \ldots,1]^T\right) \right),
\end{equation}
plays an important role as one of the building blocks of the 
proposed control laws. The vector $1_N$ is a vector of dimension $N$ with all entries equal to one. The matrix $L(a) \in \mathbb{R}^{N \times N}$ has one eigenvalue $0$ and has $N-1$ eigenvalues equal, positive and real. 

We now define the matrices  
\begin{align*}
V_1(a) & =  \left[\frac{1}{a_1}1_{N-1}, -\text{diag}\left(\left[\frac{1}{a_2}, \frac{1}{a_3},...,\frac{1}{a_N}\right]^T\right)\right], \\
V_2(a) & =  \text{diag}\left(\left[\frac{1}{a_2}, \frac{1}{a_3},...,\frac{1}{a_N}\right]^T\right) + \frac{1}{a_1}1_{N-1}1_{N-1}^{T}, \\
V_3 & =  \left[-1_{N-1}, I_{N-1}\right],
\end{align*}
and formulate the following lemma.

\begin{lem}\label{lemma:1}
$L(a) = -V_1(a)^TV_2(a)^{-1}V_3$.
\end{lem}

\quad \emph{Proof}:
A proof can be found in~\cite{thunberg2014consensus}.
\hfill $\blacksquare$ \vspace{2mm}

\begin{lem}\label{lemma:2}
Assume that $C \in \mathbb{R}^{p \times n}$ has full row rank, $P \in
\mathbb{R}^{n \times n}$ is nonsingular and
$\textnormal{ker}(C)$ is $P$-invariant, then 
$$P^TC^T(CPWP^TC^T)^{-1}CP = C^T(CWC^T)^{-1}C.$$
\end{lem}

\quad \emph{Proof}:
A proof can be found in~\cite{thunberg2014consensus}.

Let us define 
\begin{equation}
W(t,T) = \int_{t}^{T}Ce^{A(T-s)}BB^Te^{A^T(T-s)}C^Tds.
\end{equation}
The matrix $W(t,T)$ is the output controllability Gramian, and since 
the system is assumed to be output controllable, this matrix is 
nonsingular (for $t < T$).
Let us also define the related matrix 
\begin{equation}
G(t,T) = \int_{0}^{T-t}Ce^{-Ar}BB^Te^{-A^Tr}C^Tdr. \quad \: \:
\end{equation}
Beware of the difference between the transpose operator $(\cdot)^T$ and the time $T$.

The approach we use in this work relies to a large extent on the projection theorem in Hilbert spaces.
We recall the following version of the projection theorem where inner product
constraints are present.

\begin{thm}[\cite{luenberger1997optimization}]\label{thm:luen}
Let $H$ be a Hilbert space and $\{z_1, z_2, \ldots, z_N\}$ a set of linearly 
independent vectors in $H$. Among all vectors $w \in H$ satisfying 
\begin{align*}
\langle w, z_1\rangle  & = c_1, \\
\langle w, z_2\rangle  & = c_2, \\
& \vdots   \\
\langle w, z_M\rangle  & = c_M,
\end{align*}
let $z_0$ have minimum norm. Then 
$$z_0 = \sum_{i = 1}^N\beta_iz_i,$$
where the coefficients $\beta_i$ satisfy the equations
\[ \begin{array}{ccccc}
\langle z_1,z_1\rangle \beta_1 +  \langle z_2,z_1\rangle \beta_2  +  \cdots + \langle z_N,z_1\rangle \beta_N  =  c_1, & \\
\langle z_1,z_2\rangle \beta_1 +  \langle z_2,z_2\rangle \beta_2  + \cdots + \langle z_N,z_2\rangle \beta_N = c_2, & \\
 \vdots & \\
\langle z_1,z_M\rangle \beta_1 + \langle z_2,z_M\rangle \beta_2 +  \cdots + \langle z_N,z_M\rangle \beta_N  = c_M.&
\end{array}\] 

\end{thm}

In Theorem \ref{thm:luen} $\langle \cdot , \cdot \rangle $ denotes the inner product.

\section{Finite time consensus}\label{optimal1}
In this section we consider the following problem.
\begin{prob}\label{problem:1}
For any finite $T > t_0$, construct a control law $u(t)$ for the system of agents
such that the agents reach consensus in their outputs at time $T$, {i.e.}, $$y_i(T) = y_j(T) \: \: \forall i \neq j,$$ while minimizing the following cost functional 

\begin{equation}\label{eq:objective}
\int_{t_0}^{T}\sum_{i=1}^{N}a_iu_{i}^Tu_{i}dt,
\end{equation}
where $a_{i} \in \mathbb{R}^+$, $i = 1,2,...,N $.
\end{prob}

Note that the criterion in Problem \ref{problem:1} only regards
the time $T$ and does not impose any constraints on $y(t)$ when
$t > T$.
When we say that a control law $u$ for the 
system is based on only relative information, we mean 
that $$u_i = g(y_1 - y_i, \ldots, y_N - y_i), \quad \forall i,t \geq t_0$$
for some function $g$. An interesting question to answer, is under what circumstances
it is possible to construct a control law that solves
Problem~\ref{problem:1} using only relative information.
The following lemma provides a first step on the path to the answer of this question.

\begin{lem}\label{prop:proppe}
Suppose $\textnormal{ker}(C)$ is $A$-invariant and $u$ 
is based on only relative information, then
there is no locally Lipschitz continuous in state, time-invariant feedback control law $u$
that solves Problem~\ref{problem:1} and for which $g(0, \ldots, 0) = 0$. 
\end{lem}

\quad \emph{Proof of Lemma~\ref{prop:proppe}}:
Let us introduce the invertible matrix 
$$P = 
\begin{bmatrix}
C^T & C^T_{\text{ker}}
\end{bmatrix},
$$
where $C_{\text{ker}}$ has full row rank and the columns 
of $C^T_{\text{ker}}$ span $\text{ker}(C)$. Let 
us now define $\tilde{x}_i = (x_{i,1}, x_{i,2})^T$ through the following relation
$${x}_i = 
P\tilde{x}_i,
$$
for all $i$.
The dynamics for $\tilde{x}$ is given by 
$$\dot{\tilde{x}}_i = \tilde{A}\tilde{x}_i + \tilde{B}u_i, \quad y_i = \tilde{C}\tilde{x}_i,$$
where 
$$\tilde{A} = \begin{bmatrix}
A_{11} & 0  \\
A_{21} & A_{22} \\
\end{bmatrix},
\tilde{B} = \begin{bmatrix}
B_1 \\
B_2 \\
 \end{bmatrix} \text{ and }
\tilde{C} = \begin{bmatrix}
C_1 & 0
 \end{bmatrix}.
$$
The structure of $\tilde{A}$ is a consequence of the fact that $\text{ker}(C)$
is $A$-invariant. 

Suppose there is a linear time-invariant 
feedback control law $u$ that solves the Problem~\ref{problem:1}
with only relative information. 
We note that 
$$y_i = CC^Tx_{i,1}.$$
We define $y_{1j} = y_1 - y_j$ for all $j > 1$.
The control law $u$ has the following form
$$u = f(y_{12}, \ldots, y_{1N}),$$
but $y_{1j} = CC^Tx_{1j,1}$,
where $x_{1j,1} = x_{1,1} - x_{j,1}$,
so $u$ can be written as 
$$u = f(x_{12,1}, \ldots,x_{1N,1}).$$
Since $CC^T$ is invertible it follows that $y_{1j} = 0$ for all $j >1$ if and only if
$x_{1j,1} = 0$ for all $j >1$. But, by the structure of
$f$ and $\tilde{A}$ it follows that $(x_{12,1}, \ldots,x_{1N,1})^T = 0$
is an equilibrium for the dynamics of $(x_{12,1}, \ldots,x_{1N,1})^T$.
Since the right-hand side of this dynamics is locally Lipschitz continuous, 
$(x_{12,1}, \ldots,x_{1N,1})^T$ cannot have reached the point $0$ in finite time. 
This is a contradiction.
\hfill $\blacksquare$ \vspace{2mm}

We now provide
the solution to Problem \ref{problem:1}.

\begin{thm}\label{theorem:t1}
For $T < \infty$ the solution to Problem~\ref{problem:1} is 
\begin{align}
\nonumber
& u(t) = \\
\label{eq:optimal:openn1}
& -L(a) \otimes  \left(B^Te^{A^T(T-t)}C^TW(t_0,T)^{-1}Ce^{A(T-t_0)}\right)x_0\\
\nonumber
& u(x,t) =  \\
\label{eq:optimal:control1}
& -L(a) \otimes \left(B^Te^{A^T(T-t)}C^TW(t,T)^{-1}Ce^{A(T-t)}\right)x.
\end{align}
Furthermore, if $ker(C)$ is $A$-invariant, the solution to
Problem~\ref{problem:1} is  
\begin{align}
\label{eq:optimal:control}
u(y,t)  = -L(a) \otimes B^TC^TG(t,T)^{-1}y.
\end{align}
or equivalently 
\begin{align*}
u(y,t)  = B^TC^TG(t,T)^{-1}\sum_{j = 1}^N(\alpha a_jy_j - y_i), 
\end{align*}
where 
$\alpha = \left(\sum_{i=1}^Na_i \right)^{-1}.$
\end{thm}
All the control laws (\ref{eq:optimal:openn1}-\ref{eq:optimal:control}) are
equivalent (provided $\text{ker}(C)$ is $A$-invariant), but expressed in different ways.
The control law \eqref{eq:optimal:openn1} is the open loop controller and 
\eqref{eq:optimal:control1} is the closed loop version of \eqref{eq:optimal:openn1}.
The matrix $W(t,T)$ is invertible due to the assumption of
output controllability. 
We take the liberty of denoting all the controllers (\ref{eq:optimal:openn1}-\ref{eq:optimal:control})
by $u$. Provided $u$ is used during $[t_0, T)$, at the time $T$ we have that  
$$\lim_{t \uparrow T}u(x(t),t) = \lim_{t \uparrow T}u(y(t),t) = u(T).$$ 
Even though the feedback controllers in \eqref{eq:optimal:control1} and \eqref{eq:optimal:control} are bounded
and continuous for $t \in [t_0,T)$ (see the open loop version of $u$ in (\ref{eq:optimal:openn1})), computational difficulties arise as $t \rightarrow T$
when \eqref{eq:optimal:openn1} and \eqref{eq:optimal:control1} are used, since $W(T)$ is not invertible. \vspace{2mm}


\quad \emph{Proof of Theorem~\ref{theorem:t1}}:
Problem 1 is formally stated as follows
\begin{equation*}\label{eq:sing_int:3}
\text{minimize } \int_{t_0}^{T}\sum_{i=1}^{N}a_iu_{i}^Tu_idt  \quad a_{i} \in \mathbb{R}^+ \: \: i = 1,2,...,N,
\end{equation*}
when
\begin{align*}
y_i(t) & = Ce^{A(t-t_0)}x_i(t_0) + \int_{t_0}^{t}Ce^{A(t-s)}Bu_ids, \text{ for all } i, \quad
\end{align*}
and
\begin{align}
\nonumber
& \int_{t_0}^{T}Ce^{A(T-s)}B\left(u_1 - u_i\right)ds  = \\
\label{eq:constraints:initial}
&  -Ce^{A(T-t_0)}\left(x_1(t_0) - x_i(t_0)\right), 
\end{align}
for $i \in \{2,...,N\}$. 
Here we have without loss of generality assumed that the outputs of the agents at 
time $T$, $y_i(T)$ shall be equal
to the output of agent 1 at time $T$, \emph{i.e.}, $y_1(T)$.

We notice that this problem is a minimum norm problem in the Hilbert space of all functions $$f = (f_1(t), f_2(t), ..., f_{N}(t))^T: \mathbb{R} \rightarrow \mathbb{R}^{mN},$$ such that the Lebesgue integral 
\begin{equation}\label{eq:norm_L}
\int_{t_0}^{T}\sum_{i=1}^{N}a_{i}f_{i}^T(t)f_i(t)dt
\end{equation}
converges. Here $f_i:\mathbb{R} \rightarrow \mathbb{R}^m$. We denote this space $H,$ 
and the norm is given by the square root of \eqref{eq:norm_L}.

Now we continue along the lines of Theorem~\ref{thm:luen} and reformulate
the constraints \eqref{eq:constraints:initial} into inner product constraints in $H$. 
\begin{eqnarray*}\label{eq:constraints:2}
\int_{t_0}^{T}Ce^{A(T-s)}B(u_1(s) - u_i(s))ds & = & \\
\nonumber
\bigg \langle \left[\frac{Ce^{A(T-t)}B}{a_1}, 0, ...,0,-\frac{Ce^{A(T-t)}B}{a_i},0, ...,0 \right]^T, & & \\
\nonumber
\left[u_1^T, 0,...,0,u_i^T, 0, ...,0 \right]^T \bigg \rangle  = & & \\
\nonumber
-Ce^{A(T-t_0)}\left(x_1(t_0) - x_i(t_0)\right). & &
\end{eqnarray*} 
Depending on context the symbol $\langle \cdot , \cdot \rangle $ shall be interpreted as follows.
If $f$ and $g$ belongs to $H$, $\langle f,g\rangle $ denotes the inner product between these 
two elements. If $f(t)$ and $g(t)$ are matrices
of proper dimensions,
 then $\langle f,g\rangle $ is a matrix inner product where each element in the matrix is 
an inner product between a column in $f$ and a column in $g$.

To simplify the notation we define 
$$p_i = \left[\frac{Ce^{A(T-t)}B}{a_1}, 0,...,0,-\frac{Ce^{A(T-t)}B}{a_i},0, ...,0 \right]^T.$$ 
Since we have a minimum norm problem and all the columns of all the $p_i$ are independent, 
by Theorem~\ref{thm:luen} we get that the optimal controller $u(t)$ is given by
\begin{equation}\label{eq:100}
u(t) = [p_2, ...,p_N]\beta,
\end{equation}
where $\beta$ is the solution to
\begin{equation}
Q\beta = V_3 \otimes Ce^{A(T-t_0)}x_0,
\end{equation}
where
\begin{align}
\label{eq:101}
Q = \begin{bmatrix}
\langle p_2,p_2\rangle  & \langle p_3,p_2\rangle  & \cdots &\langle p_N,p_2\rangle  \\
\langle p_2,p_3\rangle  & \langle p_3,p_3\rangle  & \cdots &\langle p_N,p_3\rangle  \\
\vdots & \vdots & \ddots & \vdots & \\
\langle p_2,p_N\rangle  & \langle p_3,p_N\rangle  & \cdots &\langle p_N,p_N\rangle  
 \end{bmatrix}.
\end{align}
From (\ref{eq:100}-\ref{eq:101}) we get that $\beta = Q^{-1}V_3 \otimes Ce^{A(T-t_0)}x_0$ and $u = [p_2, ...,p_N] Q^{-1}V_3 \otimes Ce^{A(T-t_0)}x_0$.
Now we have that $[p_2, ...,p_N] = V_1(a)^T \otimes (B^Te^{A^T(T-t)}C^T)$. 

Since
\begin{equation*}
\langle p_i,p_j\rangle  =
\begin{cases}
\frac{1}{a_1} W(t_0,T) \quad & \text{if} \: i \neq j, \\
\left(\frac{1}{a_1} + \frac{1}{a_{i+1}}\right) W(t_0,T) \quad & \text{if} \: i = j,
\end{cases}
\end{equation*}
where $W(t_0,T) = \int_{t_0}^{T}Ce^{A(T-s)}BB^Te^{A^T(T-s)}C^Tds$, \\ we have that 
\begin{equation*}
Q = V_2(a) \otimes W(t_0,T).
\end{equation*}
\begin{align*}
\nonumber
u(t)  & =  \left( V_1(a)^T  \otimes  (B^Te^{A^T(T-t)}C^T) \right)\cdot \\
& \left( V_2(a)^{-1}  \otimes
 W(t_0,T)^{-1} \right) \left( V_3 \otimes Ce^{A(T-t_0)}\right) x_0 \\
 & =  \left(V_1(a)^TV_2(a)^{-1}V_3\right) \otimes \\
 & \quad \:\left(B^Te^{A^T(T-t)}C^TW(t_0,T)^{-1}Ce^{A(T-t_0)}\right)x_0   \\
 & = -L(a) \:\otimes \\
 &  \quad \:\left( B^Te^{A^T(T-t)}C^TW(t_0,T)^{-1}Ce^{A(T-t_0)}\right)x_0.
\end{align*}
By Bellman's Principle we get  that
\begin{align*}
& u(x,t)  =  \\
&  -L(a) \otimes \left(B^Te^{A^T(T-t)}C^TW(t,T)^{-1}Ce^{A(T-t)}\right)x.
\end{align*}
Now, provided $\text{ker}(C)$ is $A$-invariant, we can use 
Lemma~\ref{lemma:2} to obtain
\begin{align}
\label{eq:u:y}
u(y,t)  = -L(a) \otimes B^TC^TG(t,T)^{-1}y
\end{align}
(see \cite{thunberg2014consensus} for details).
\hfill $\blacksquare$ \vspace{2mm}


\begin{cor}
The controllers (\ref{eq:optimal:control1}-\ref{eq:optimal:control}) use only relative information, {i.e.},~differences of the states (outputs)
of the agents.
\end{cor}

\quad \emph{Proof}:
Straight forward by using the structure of the matrix $L(a)$.
\hfill $\blacksquare$ \vspace{2mm}

Let us define $y_c=\frac{1}{\sum_{i=1}^Na_i}\sum_{i=1}^Na_iy_i,$
  \:  and $\bar{y}_c  = (y_c, \ldots, y_c)^T \in \mathbb{R}^{pN}$.

\begin{lem}\label{theorem:meeting_point}
Suppose that $A$ has not full rank and $x_i(0) = x_{i0} \in \text{ker}(A)$ for all $ i=1,...,N$, then the consensus point 
for the system of agents using the controller \eqref{eq:optimal:control1} or \eqref{eq:optimal:control} is $y_c(0)$.
\end{lem}

\quad \emph{Proof}:
We have that 
\begin{align*}
y(T) = I_N \otimes Ce^{A(T-t_0)}x_0 + \int_{t_0}^T\left(I_N \otimes Ce^{A(T-t)}B\right)\cdot & \\
\left(-L(a) \otimes B^Te^{A^T(T-t)}C^TW(t_0,T)^{-1}Ce^{A(T-t_0)}\right)x_0dt & \\
 = y_0 + \int_{t_0}^T (-L(a)\otimes & \\
 Ce^{A(T-t)}BB^Te^{A^T(T-t)}C^TW(t_0,T)^{-1} )y_0dt = \bar{y}_c(t_0). &
\end{align*} 
\hfill $\blacksquare$ \vspace{2mm}

\begin{thm}
For non-homogeneous output controllable $(A_i,B_i,C_i),~i=1,\cdots,N$, the optimal
controller can be derived in a similar way and has the following form:
$$u_i(t) = B_i^Te^{A_i^Tt}C_i^TW_i(0,T)^{-1}(\alpha^*-
C_ie^{A_it}x_i(0)),~i\le N,$$
where 
$$ W_i(0,T) = \int_0^TC_ie^{A_i(T-s)}B_iB_i^Te^{A_i(T-s)}C_i^Tds,$$

\begin{align*}
\alpha^* = &   \left(\sum_{j = 1}^Na_jW_j(0,T)^{-1}\right)^{-1} \cdot \\
& \sum_{i
  = 1}^N{a_i}W_i(0,T)^{-1}C_ie^{A_it}x_i(0).
\end{align*}
\end{thm}

\quad \emph{Proof}:\\
\textbf{Step 1:} \\
Solve the problem 
\begin{equation*}
\begin{cases}
& \min\limits_{u_i}\int\limits_0^T\frac{1}{2}u_i(t)u_i(t)dt \\
& \text{s.t.}~\dot{x} = A_ix + B_iu, \\
& \quad \: \: \: y_i(T) = \alpha.
\end{cases}
\end{equation*} \\
Using any standard approach, one obtains that 
$$u_i^*(t,\alpha) = B_i^Te^{A_i^Tt}C_i^TW_i(0,T)^{-1}(\alpha - C_ie^{A_it}x_i(0))$$
is the solution to this problem provided ($A_i$,$B_i$,$C_i$) is output controllable. 

\textbf{Step 2:} \\
We want to find the $\alpha^*$ that minimizes 
$$\int_0^T\sum_{i = 1}^Na_iu^*_i(t,\alpha)^Tu_i^*(t,\alpha)dt.$$
Now,
\begin{align*}
& \min_{\alpha}\int_0^T\sum_{i = 1}^Na_iu^*_i(t,\alpha)^Tu_i^*(t,\alpha)dt \\
 = & 
\min_{\alpha}\sum_{i = 1}^Na_i(\alpha -C_ie^{A_it}x_i(0))^TW_i(0,T)^{-1}\cdot \\
& (\alpha -C_ie^{A_it}x_i(0)).
\end{align*}
This gives that 
$$\sum_{i = 1}^N{a_i}W_i(0,T)^{-1}(\alpha^* - Ce^{A_it}x_i(0)) = 0$$
or 
\begin{align*}
\alpha^* = &  \left(\sum_{j = 1}^Na_jW_j(0,T)^{-1}\right)^{-1}\cdot \\ 
& \sum_{i = 1}^N{\alpha_i}W_i(0,T)^{-1}C_ie^{A_it}x_i(0).
\end{align*}
Thus the optimal controller $u = [u_1^T, u_2^T, \ldots, u_N^T]$ is given by 
\begin{align*}
u_i(t) =&  B_i^Te^{A_i^Tt}C_i^TW_i(0,T)^{-1}\cdot \\
& (\frac{1}{\sum_{j = 1}^Na_j}\sum_{i = 1}^N{a_i}C_ie^{A_it}x_i(0) - Ce^{A_it}x_i(0))
\end{align*}
for all $i$.
\hfill $\blacksquare$ \vspace{2mm}

\section{Extension to the asymptotic consensus problem}\label{optimal2}
We now examine the asymptotic case,
\emph{i.e.}, we want the system to asymptotically reach consensus while
minimizing the cost functional. The problem is formally 
stated as follows.

\begin{prob}\label{problem:2}
Construct a control law $u(t)$ for the system of agents
such that the agents asymptotically reach consensus in the outputs, \emph{i.e.},
$$\lim_{t \rightarrow \infty} (y_i(t) - y_j(t)) = 0 \: \: \text{ for all } i,j \text{ such that } i \neq j,$$ while minimizing the following cost functional 

\begin{equation}
\int_{t_0}^{\infty}\sum_{i=1}^{N}a_iu_{i}^T(t)u_{i}(t)dt
\end{equation}
where $a_{i} \in \mathbb{R}^+$ for $i = 1,2,...,N $.
\end{prob}
In order to solve Problem~\ref{problem:2},
we start by defining the matrix
$$P(t,T)=e^{A^T(T-t)}C^TW(t,T)^{-1}Ce^{A(T-t)}$$ which satisfies
the following differential Riccati equation
\begin{equation}\label{ricc1}
\frac{dP}{dt}=-A^TP-PA+PBB^TP.
\end{equation}
The matrix $P(t,T)$ is an essential part of the control laws
that were presented in the last section, and here we see that
this matrix is provided as the solution to a differential matrix Riccati equation.
It is well known that \eqref{ricc1} has a positive semidefinite limit $P_0$ as $T-t\rightarrow\infty$ if 
$(A,B)$ is stabilizable and $A$ does not have any
eigenvalue on the imaginary axis. In order to see this
we consider the following problem
\begin{equation}
\label{eq:other}
\begin{aligned}
\text{min} &~\int_{0}^{\infty}\Vert u\Vert^2 dt\\
s.t. &~\dot x=Ax+Bu.
\end{aligned}
\end{equation}
If $(A,B)$ is stabilizable and $A$ does not have any
eigenvalue on the imaginary axis,
$$u=-B^TP_0x$$
is the optimal control law that solves \eqref{eq:other},
where $P_0$ is the positive semi-definite solution to
$$-A^TP_0-P_0A +P_0BB^TP_0=0.$$
This Algebraic Riccati equation is obtained by letting 
the left-hand side of \eqref{ricc1} be equal to zero.

The problem \eqref{eq:other} is not a consensus problem,
and the question is, besides the fact the 
same matrix $P_0$ is used in the optimal control law, how it is related to our consensus problem. It turns out that the 
control law, besides being the solution of the consensus problem, is also the solution of $N$ problems on the form \eqref{eq:other}.
In order to show this we introduce
$$x_c=\frac{1}{\sum_{i=1}^Na_i}\sum_{i=1}^Na_ix_i, \text{ and } \delta_i=x_i-x_c.$$
The dynamics of $x_c$ and $\delta_i$ are given by
$$\dot x_c=Ax_c \quad \text{ and } \quad \dot \delta_i=A\delta_i+Bu_i.$$
Now each control law $u_i(t)$ contained in the vector
$$u(t)= -L(a) \otimes (B^TP_0)x,$$
can be written as 
$$u_i = B^TP_0\delta_i$$
where $u_i$ solves the problem 
\begin{equation}
\label{eq:other2}
\begin{aligned}
\text{min} &~\int_{0}^{\infty}\Vert u\Vert^2 dt\\
s.t. &~\dot \delta_i=A\delta_i+Bu_i.
\end{aligned}
\end{equation}
Provided $\text{ker}(C)$ is $A$-invariant, it can be shown that 
\begin{align*}
P(t,T)& =e^{A^T(T-t)}C^TW(t,T)^{-1}Ce^{A(T-t)} \\
& = C^TG(t,T)^{-1}C,
\end{align*}
which implies that 
$$(CC^T)^{-1}CP(t,T)C^T(CC^T)^{-1} = G(t,T)^{-1}.$$
Now, as $T-t \rightarrow \infty$, it holds that
$P(t,T) \rightarrow P_0$. This means that, as $T-t \rightarrow \infty$, 
\begin{align*}
(CC^T)^{-1}CP(t,T)C^T(CC^T)^{-1} \rightarrow \\ (CC^T)^{-1}CP_0C^T(CC^T)^{-1}.
\end{align*}

Let $G_0 = (CC^T)^{-1}CP_0C^T(CC^T)^{-1}$.
Then for the asymptotic consensus problem,
controller \eqref{eq:u:y} becomes
\begin{equation*}\label{nisse:47}
u = -L(a) \otimes (B^TC^TG_0)y.
\end{equation*}

\begin{prop}
If $A$ is stabilizable with 
no eigenvalues on the imaginary axis. Then
$P_0$ exists, is positive semidefinite and the optimal control law that solves Problem~\ref{problem:2} is
$$u= -L(a) \otimes (B^TP_0)x.$$
Furthermore, if $\text{ker}(C)$ is also A-invariant then 
$$u= -L(a) \otimes (B^TC^TG_0)y.$$
\end{prop}

When only the output $y_i=Cx_i$ is available for control action and $\text{ker}(C)$ is not necessarily $A$-invariant, an
observer can be designed.
$$\dot{\hat\delta}_i=(A-BB^TP_0)\hat\delta_i-QC^T\left(\frac{1}{\sum_{i=1}^Na_i}\sum_{j=1}^Na_j(y_i-y_j)-C\hat\delta_i\right).$$

Under the assumption that $(A,C)$ is detectable and $A$ does not have any
eigenvalue on the imaginary axis we have that 
\begin{equation*}
\dot{\hat \delta}=A\hat \delta -BB^TP_0\hat \delta -QC^T((y_i - y_c) -C\hat \delta),
\end{equation*}
where $Q\le 0$ satisfies
$$AQ+QA^T = -QC^TCQ.$$
We summarize these results in the following proposition.

\begin{prop}
Suppose $(A,B)$ is stabilizable and $(A,C)$ is detectable, and $A$ has
no eigenvalue on the imaginary axis. Then, if the following dynamic output control law
is used,
\begin{align*}
u_i & =  -B^TP_0\hat\delta_i,\\
& \dot{\hat\delta}_i = (A-BB^TP_0)\hat\delta_i - \\
\nonumber
&   QC^T\left(\frac{1}{\sum_{i=1}^Na_i}\sum_{j=1}^Na_j(y_i-y_j)-C\hat\delta_i\right),
\end{align*}
the system reaches asymptotic consensus in the outputs.
\end{prop}

\bibliographystyle{agsm}        
\bibliography{ref}

\end{document}